\documentclass[12pt]{amsart}

\usepackage{fullpage} 
\usepackage{parskip} 
\usepackage{tikz-cd} 
\usepackage{amsmath, amssymb, amsthm}
\usepackage{esint}
\usepackage{hyperref}
\usepackage[utf8]{inputenc}
\usepackage[english]{babel}
\usepackage{amsrefs}

\usepackage[title]{appendix}

\usepackage{csquotes}

\usepackage{geometry}
\usepackage{fancyhdr}
\usepackage{lipsum}

\numberwithin{equation}{section}

\theoremstyle{definition}
\newtheorem{thm}[equation]{Theorem} 
\theoremstyle{definition}

\newtheorem{lemma}[equation]{Lemma}

\newcommand{\ovl}[1]{\overline{#1}}
\newcommand{\pair}[1]{\left\langle #1\right\rangle}
\newcommand{\bb}[1]{\mathbb{#1}}

\newcommand\numberthis{\addtocounter{equation}{1}\tag{\theequation}}

\newcommand{\tr}{\text{tr}}
\newcommand{\D}{\mathcal{L}}
\newcommand{\Hess}{\text{Hess}}
\newcommand{\n}{\textbf{n}}

\title{Convexity of $\lambda$-hypersurfaces}
\author{Tang-Kai Lee}
\address{MIT, Dept. of Math., 77 Massachusetts Avenue, Cambridge, MA 02139-4307}
\email{tangkai@mit.edu}

\date{\today}

\begin{document}
	\maketitle
	\begin{abstract}
	We prove that any $n$-dimensional closed mean convex $\lambda$-hypersurface is convex if $\lambda\le 0.$ This generalizes Guang's work on $2$-dimensional strictly mean convex $\lambda$-hypersurfaces. As a corollary, we obtain a gap theorem for closed $\lambda$-hypersurfaces with $\lambda\le 0.$
	\end{abstract}
	\section{{\bf Introduction}}
A hypersurface $M^n$ in $\bb R^{n+1}$ is called a $\lambda$-\textit{hypersurface} if it satisfies
\begin{equation}\label{lambda}
H - \frac{\pair{x,\n}}{2} = \lambda
\end{equation}
where $H$ is the mean curvature, $\n$ is the outer unit normal of $M,$ $x$ is the position vector, and $\lambda$ is a constant. This equation arises in the study of isoperimetric problems in weighted (Gaussian) Euclidean spaces (c.f. \cite{MR15}), which is a long-standing topic studied in various fields in science (\cite{L94}, \cite{B01}, \cite{B03}, etc.). Recently, Cheng and Wei \cite{CW18} defined a weighted volume functional, and showed that the critical points of the functional under some weighted volume-preserving variations are exactly $\lambda$-hypersurfaces.

When $\lambda=0,$ $\lambda$-hypersurfaces are exactly \textit{self-shrinkers}. Self-shrinkers play an important role in the study of mean curvature flow (MCF), since White \cite{W97} and Ilmanen \cite{I95} showed that self-shrinkers arise as the tangent flows of MCF based on Huisken's monotonicity formula \cite{H90} and Brakke's compactness theorem \cite{B78}. Many classification results of self-shrinkers were proposed. Abresch and Langer \cite{AL86} showed that the only $1$-dimensional closed embedded self-shrinker is the circle $S^1.$ Huisken \cite{H90} later dealt with the higher-dimensional cases, proving that any closed, embedded, and mean convex (which means $H\ge 0$) $n$-dimensional self-shrinkers are exactly spheres $S^n.$ For the non-compact situation, Huisken \cite{H93} proved that all smooth, embedded, and mean convex self-shrinkers with polynomial volume growth and bounded second fundamental form are generalized cylinders $S^k\times \bb R^{n-k}.$ This result was later improved by Colding and Minicozzi \cite{CM12}, in which they removed the condition of bounded second fundamental form in Huisken's classification. 

When $\lambda\neq 0,$ there are relatively few and incomplete classification results so far. In \cite{CW18}, Cheng and Wei characterized compact $\lambda$-hypersurfaces with $H-\lambda\ge 0$ and some curvature conditions (c.f. theorem \ref{CW}). Inspired by \cite{SX20}, Guang \cite{G21} showed that any strictly mean convex (which means $H>0$) $2$-dimensional $\lambda$-hypersurfaces are in fact convex if $\lambda\le 0.$ The main goal of this paper is to generalize Guang's result to higher-dimensional mean convex $\lambda$-hypersurfaces. 

\begin{thm}\label{main}
	Let $M^n$ be a smooth, closed, and embedded $\lambda$-hypersurface in $\bb R^{n+1}$ with $\lambda\le 0.$ If $M$ is mean convex, then it is convex.
\end{thm}

This theorem is a generalization of Guang's result in \cite{G21}. Guang used the explicit expressions for the derivatives of the principal curvatures at the non-umbilical points of a surface, which were first derived in \cite{HIMW19}. We follow the same spirit to derive a differential inequality for the sum of a part of principal curvatures at the points where there is a gap among some principal curvatures (c.f. lemma \ref{LS}). Though we could not derive similar explicit expressions, it turns out that the information we derive is sufficient to obtain the higher-dimensional generalization. Also, we use the maximum principle to weaken the assumption of strict mean convexity in \cite{G21}, which Huisken \cite{H90} also applied when classifying closed mean convex self-shrinkers.

A natural further question is whether Huisken-type classification also holds for $\lambda$-hypersurfaces. That is, could we classify all $\lambda$-hypersurfaces given some curvature conditions, like mean-convexity? In the curve case, Guang \cite{G18} proved that any smooth embedded $1$-dimensional $\lambda$-hypersurface (or $\lambda$-curve) is either a straight line or a circle if $\lambda\ge 0,$ which generalized Abresch and Langer's result. For the higher dimensional case, Heilman \cite{H17} proved that convex $n$-dimensional $\lambda$-hypersurfaces are generalized cylinders if $\lambda\ge 0.$ However, when $\lambda<0,$ Chang \cite{C17} showed that for certain $\lambda<0,$ there are some closed embedded mean convex $\lambda$-curves other than circles. Thus we could not expect Huisken-type results to hold for general $\lambda\in\bb R.$ We hope that theorem \ref{main} will shed some light on the higher-dimensional case when $\lambda\le 0.$ In particular, using the curvature condition discovered in \cite{CW18}, we can prove the following gap theorem for mean convex $\lambda$-hypersurfaces when $\lambda\le 0.$

\begin{thm}\label{gap}
	Let $M^n$ be a smooth, closed, and embedded $\lambda$-hypersurface in $\bb R^{n+1}.$ If $\lambda\le 0$ and the mean curvature of $M$ satisfies
	$$0\le H\le \frac{\sqrt{\lambda^2+2}+\lambda}{2},$$
	then $M$ is a round sphere.
\end{thm}

We remark that if we assume $M$ is a convex $\lambda$-hypersurface, then the result of theorem \ref{gap} could also be derived from the gap theorem proven by Guang \cite{G18}. What's new here is that we only need to assume $H\ge 0,$ and then by theorem \ref{main}, we can get the convexity. 

The organization of this paper is as follows. In section \ref{2}, we will introduce the Simons-type identities for $\lambda$-hypersurfaces, which was given by Guang in \cite{G18}. In section \ref{3}, we derive a differential inequality for the sum of a part of principal curvatures at the points where there is a gap among some principal curvatures. In section \ref{4}, we use the identities and the inequality in the preceding sections to prove the main theorem \ref{main}. In section \ref{5}, we prove a gap theorem \ref{gap} by applying Cheng and Wei's theorem.

\subsection*{\bf Acknowledgement}
The author is grateful to Prof. Bill Minicozzi for his helpful and inspiring comments. He also appreciates Kai-Hsiang Wang's indications on some deficiencies in an earlier draft, and Qiang Guang's generosity on sharing some useful references. This work was completed when the author visited the National Center for Theoretical Science (NCTS) in Taiwan, and the author is also grateful for much helpful discussion with people in NCTS.

\section{\bf Simons-type Identities}\label{2}
On a hypersurface $M$ in $\bb R^{n+1},$ we consider the drift Laplacian
$$\D := \Delta-\frac 12\nabla_{x^T}(\cdot)$$ 
and also the following linear operator
$$L := \D + |A|^2 + \frac 12 
= \Delta - \frac 12\nabla_{x^T}(\cdot) + |A|^2 + \frac 12$$
where $\Delta$ and $A$ denote the Laplacian operator and the second fundamental form of $M,$ and $x^T$ is the tangential component (with respect to $M$) of the position vector $x.$ These operators were introduced by Colding and Minicozzi to study the stability of self-shrinkers. In fact, the operator $L$ appears in the second variation formula of the $F$-functional (c.f. \cite{CM12}).

Guang \cite{G18} established the following Simons-type identities. These identities will play a crucial role in the proof of the main theorem \ref{main}. We remark that these kinds of identities have been developed in \cite{CM12} and \cite{CM15} for self-shrinkers. For completeness, we include the proof in \cite{G18} here.
\begin{lemma}[\cite{G18}]\label{LALH} If $M$ is a $\lambda$-hypersurface in $\bb R^{n+1},$ then
	\begin{equation}\label{LA}
	LA = A-\lambda A^2
	\end{equation}
	and in particular, taking the trace of \eqref{LA} gives
	\begin{equation}\label{LH}
	LH = H + \lambda|A|^2.
	\end{equation}
\end{lemma}

\noindent{\bf (Proof.)} For any fixed $p\in M,$ take a local orthonormal frame $\{e_i\}_{i=1,\cdots,n}$ such that $\nabla^M_{e_i}e_j=0$ for all $i$ and $j,$ where $\nabla^M$ is the Riemannian connection of $M.$ Thus we can write $\nabla_{e_i}e_j=a_{ij}n$ where $a_{ij}$ is the component of the second fundamental form $A.$ As a result,
\begin{align*}
\Hess_{\pair{x,\n}}(e_i,e_j)
= \nabla_{e_j}\nabla_{e_i}\pair{x,\n}
& = \nabla_{e_j}\sum_{l=1}^n \pair{x,-a_{il}e_l}\\
& = -a_{ij} + \sum_{l=1}^n \left(-a_{il,j}\pair{x,e_l} - a_{il}\pair{x,a_{jl}\n} \right)\\
& = -A(e_i,e_j) - (\nabla_{x^T}A)(e_i,e_j) - \pair{x,\n}A^2(e_i,e_j)
\end{align*}
where $a_{il,j}$ is the component of $\nabla A,$ and we use the Codazzi equation $a_{il,j} = a_{ij,l}.$ In conclusion, we derive
\begin{equation}\label{Hij}
\Hess_{\pair{x,\n}} = -A - \nabla_{x^T}A - \pair{x,\n}A^2.
\end{equation}
Plug this into the Simons identity
\begin{equation}\label{Simons}
\Delta A = -|A|^2 A - HA^2 - \Hess_H
\end{equation}
which holds for any hypersurface in $\bb R^n$ (c.f. the formula (2.14) in \cite{CM11}), and we get
\begin{align*}
LA
& = \Delta A - \frac 12\nabla_{x^T}(A) + |A|^2A + \frac 12 A\\
& = A - \left(H - \frac {\pair{x,\n}}{2}\right)A^2\\
& = A - \lambda A^2
\end{align*}
based on the $\lambda$-hypersurface equation \eqref{lambda}. \eqref{LH} follows directly after taking the trace since $\text{tr} A = - H.$\qed

\section{\bf Estimates of Principal Curvatures}\label{3}
In this section, we let $M$ be a smooth mean convex hypersurface in $\bb R^{n+1}.$ Besides, we will write $k_1\le\cdots\le k_n$ to be the principal curvatures of $M$ in the ascending order. For $l\ge 1,$ consider 
\begin{equation*}
S_l:=\sum_{m=l+1}^nk_m,
\end{equation*}
which is the sum of the largest $n-l$ principal curvatures. In general, $S_l$ is just a continuous function on $M.$ However, if $k_l<k_{l+1}$ at a point $p\in M,$ the inverse function theorem will imply $k_1+\cdots+k_l$ and thus $S_l$ are both differentiable near $p.$ At such a point, we establish the following differential inequality for $S_l.$
\begin{lemma}\label{LS}
	Suppose $k_l<k_{l+1}$ for some $l\ge 1$ at a point $p\in M.$ Then at $p,$ we have
	\begin{equation}\label{LS>=}
	\D S_l \ge \frac {S_l}2 - |A|^2 S_l +\lambda\sum_{m=l+1}^nk_m^2.
	\end{equation}
\end{lemma}
\noindent{\bf(Proof.)} We only need to consider those points near which we could take a principal frame $\{v_1,\cdots,v_n\}$ such that
\begin{equation}\label{ki}
k_i=-a_{ii}:=-A(v_i,v_i)
\end{equation}
and 
\begin{equation}\label{aperp}
a_{ij}:=A(v_i,v_j)=0\text{ for }1\le i\neq j\le n.
\end{equation}
Such points form a dense and open set in $M$ (c.f. \cite{S75}), so after proving \eqref{LS>=} at these points, it follows that \eqref{LS>=} holds for all $p\in\{k_l<k_{l+1}\}$ by continuity.

Now assume $v_1,\cdots, v_n$ form a principal frame near $p.$ For any fixed $i,$ since $\pair{v_i,v_i}=1,$ we have
\begin{equation}\label{perp}
\pair{\nabla_v v_i,v_i}=\frac 12\nabla_v\pair{v_i,v_i}=0
\end{equation}
for any local vector field $v.$ Hence we can write
\begin{equation}\label{cimj}
\nabla_{v_i}v_m=\sum_{j\neq m}c_{i}^{mj}v_j
\end{equation}
for some smooth functions $c_{i}^{mj}$ near $p.$ Based on \eqref{aperp} and \eqref{cimj}, near the point $p,$ we have
\begin{equation*}
\nabla_{v_i}k_m
= -\nabla_{v_i}(A(v_m,v_m))
= -(\nabla_{v_i}A)(v_m,v_m) + 2A(\nabla_{v_i} v_m,v_m)
= -(\nabla_{v_i}A)(v_m,v_m),
\end{equation*}
so
\begin{align*}
\Delta k_m
 = -\sum_{i=1}^n\nabla_{v_i}\nabla_{v_i}(A(v_m,v_m))
& = -\sum_{i=1}^n\nabla_{v_i}((\nabla_{v_i}A)(v_m,v_m))\\
& = - (\Delta A)(v_m,v_m) - 2\sum_{i=1}^n(\nabla_{v_i}A)(\nabla_{v_i}v_m,v_m)\\
& = - (\Delta A)(v_m,v_m) - 2\sum_{i=1}^n\sum_{j\neq m}c_i^{mj}a_{jm,i}\numberthis\label{deltak_m}
\end{align*}
by \eqref{cimj}. To calculate the term involving the derivative of the second fundamental form, notice that for $j\neq m,$ $A(v_j,v_m)=0,$ based on which we have
\begin{align*}
0
& = \nabla_{v_i}(A(v_j,v_m))\\
& = a_{jm,i} + A(\nabla_{v_i}v_j,v_m) + A(v_j,\nabla_{v_i}v_m)\\
& = a_{jm,i} + A\left(\sum_{l\neq j}c_i^{jl}v_l,v_m \right) + A\left(v_j, \sum_{l\neq m}c_i^{ml}v_l \right)\\
& = a_{jm,i} - c_i^{jm}k_m - c_i^{mj}k_j\numberthis\label{half}
\end{align*}
where we use the decomposition \eqref{cimj} and the relations \eqref{ki} and \eqref{aperp}. To get a more precise form, observe that the orthogonality condition $\pair{v_j,v_m}=0$ implies
\begin{equation}\label{ccommutator}
0
= \nabla_{v_i}\pair{v_j,v_m}
= \pair{\nabla_{v_i}v_j,v_m} + \pair{v_j,\nabla_{v_i}v_m}
= c_i^{jm} + c_i^{mj}
\end{equation}
due to the orthonormality and \eqref{cimj}. Putting \eqref{ccommutator} back into \eqref{half}, we obtain
$$0 = a_{jm,i} + c_i^{mj}(k_m-k_j),$$
with which we could simplify \eqref{deltak_m} as
\begin{equation}\label{dkm}
\Delta k_m
= - (\Delta A)(v_m,v_m) + 2\sum_{i=1}^n\sum_{j\neq m} (c_i^{mj})^2(k_m-k_j).
\end{equation}
Now we apply the Simons-type identity \eqref{LA}, which gives
\begin{align*}
(\Delta A)(v_m,v_m) 
& = \frac 12 (\nabla_{x^T}A) (v_m,v_m)
+ \frac 12A(v_m,v_m) 
- |A|^2A(v_m,v_m) 
- \lambda A^2(v_m,v_m).\\
& = - \frac 12 \nabla_{x^T} k_m
- \frac 12 k_m
+ |A|^2 k_m 
- \lambda k_m^2.
\end{align*}
Combining this with \eqref{dkm}, we derive
\begin{equation*}
\Delta k_m
= \frac 12 \nabla_{x^T} k_m
+ \frac 12 k_m
- |A|^2 k_m 
+ \lambda k_m^2
+ 2\sum_{i=1}^n\sum_{j\neq m} (c_i^{mj})^2(k_m-k_j).
\end{equation*}
As a result,
\begin{equation*}
\D k_m
= \Delta k_m - \frac 12\nabla_{x^T} k_m
= \frac 12 k_m
- |A|^2 k_m 
+ \lambda k_m^2
+ 2\sum_{i=1}^n\sum_{j\neq m} (c_i^{mj})^2(k_m-k_j).
\end{equation*}
Therefore, summing over $m$ from $l+1$ to $n$ leads to
\begin{align*}
\D S_l
= \sum_{m=l+1}^n\D k_m
& = \frac 12 \sum_{m=l+1}^n k_m 
- |A|^2 \sum_{m=l+1}^n k_m
+ \lambda \sum_{m=l+1}^n k_m^2
+ 2\sum_{m=l+1}^n\sum_{i=1}^n\sum_{j\neq m} (c_i^{mj})^2(k_m-k_j)\\
& = \frac 12 S_l - |A|^2S_l + \lambda \sum_{m=l+1}^n k_m^2
+ 2 \sum_{i=1}^n \sum_{m=l+1}^n\sum_{j=1}^l(c_i^{mj})^2(k_m-k_j) 
\end{align*}
where some of the terms in the large sum get cancelled when $m$ and $j$ are switched since \eqref{ccommutator} implies
$$(c_i^{jm})^2 = (c_i^{mj})^2$$
for all $j\neq m.$ Then the inequality \eqref{LS>=} follows since by our convention, $k_m-k_j\ge 0$ for all $m>l\ge j.$\qed

\section{\bf Proof of the Main Theorem}\label{4}
We are in a position to prove the main theorem \ref{main} using lemma \ref{LALH} and \ref{LS}. We state the main theorem here again.
\begin{thm}
	Let $M^n$ be a smooth, closed, and embedded $\lambda$-hypersurface in $\bb R^{n+1}$ with $\lambda\le 0.$ If $M$ is mean convex, then it is convex.
\end{thm}
\noindent{\bf(Proof.)} The case with $\lambda=0$ directly follows from the classification of closed mean convex self-shrinkers, so we may assume $M$ is a mean convex $\lambda$-hypersurface with $\lambda<0.$

First we show that $M$ is strictly convex. In fact, \eqref{LH} implies
\begin{align*}
\Delta H - \frac 12 \nabla_{x^T} H + \left(|A|^2 - \frac 12 \right)H
= \lambda |A|^2
\le 0.
\end{align*}
Therefore, if $H$ vanished at some points, the maximum principle would imply that $H\equiv 0.$ Thus $M$ would be planar, contradicting the assumption. Consequently we verify that $M$ is strictly convex. That is, $H>0$ on $M.$ In particular, $S_l>0$ on $M$ for all $l\ge 1.$

Next, we will prove the conclusion of the theorem by a contradiction argument. That is, assume there existed $\ovl p\in M$ such that $k_1(\ovl p)<0.$ Then 
$$\frac H{S_1}=1+\frac {k_1}{S_1}$$
would attain its minimum at such point, say at $p.$ We can find $l\ge 1$ such that at this point $p,$
$$k_1=\cdots=k_l<k_{l+1}.$$
We claim that at $p,$ the function $\frac H{S_l}$ also attains its minimum. Otherwise, if 
$\frac{H(q)}{S_l(q)}<\frac{H(p)}{S_l(p)}$	
for some $q\neq p,$ which means
$$\frac{\sum_{m=1}^l k_m(q)}{H(q)-\sum_{m=1}^l k_m(q)} < \frac{\sum_{m=1}^l k_m(p)}{H(p)-\sum_{m=1}^l k_m(p)},$$
then after expanding the terms, we get
$$H(p)\sum_{m=1}^l k_m(q)
<H(q)\sum_{m=1}^l k_m(p)
= H(q)\cdot lk_1(p).$$
This particularly implies
$$H(p)k_1(q)
\le H(p)\cdot \frac 1l\sum_{m=1}^l k_m(q)
< H(q)k_1(p),
$$
which then results in 
$$\frac {H(q)}{S_1(q)} 
= 1 + \frac{k_1(q)}{S_1(q)}
< 1 + \frac{k_1(p)}{S_1(p)}
= \frac {H(p)}{S_1(p)},
$$
contradicting the minimality of $\frac H{S_1}$ at $p.$ Thus we prove that $\frac H{S_l}$ attains its minimum at $p.$ Consequently, we have
\begin{equation}\label{41}
\D \left(\frac H{S_l} \right) \ge 0\text{ and }\nabla \left(\frac H{S_l} \right)=0
\end{equation}
at $p,$ where $S_l$ is differentiable at $p$ since $k_l(p)<k_{l+1}(p).$ Note that \eqref{LH} implies
$$\D H = \frac H2 - |A|^2 H + \lambda|A|^2.$$
Combining this with lemma \ref{LS}, we obtain that at $p,$
\begin{align*}
\D \left(\frac H{S_l} \right)
& = \frac{S_l\D H - H\D S_l}{S_l^2} - 2\pair{\nabla\left(\frac H{S_l} \right), \frac{\nabla S_l}{S_l}}\\
& \le \frac{1}{S_l}\left(\frac H2 - |A|^2 H + \lambda|A|^2 \right) - \frac H{S_l^2}\left(\frac {S_l}2 - |A|^2 S_l +\lambda\sum_{m=l+1}^nk_m^2 \right)\\
& = \frac \lambda {S_l} \left(|A|^2 - \frac H{S_l} \sum_{m=l+1}^nk_m^2 \right)\\
& = \frac{\lambda}{S_l} \left(\sum_{i=1}^n k_i^2 -\left(1+\frac {\sum_{j=1}^lk_j}{S_l} \right) \left(\sum_{m=l+1}^n k_m^2 \right)\right)\\
& = \frac{\lambda}{S_l}\left(\sum_{i=1}^lk_i^2 - \frac {\sum_{j=1}^lk_j}{S_l}\left(\sum_{m=l+1}^n k_m^2 \right)\right)\\
& = \frac{l\lambda k_1}{S_l}\left(k_1 -\frac 1{S_l}\left(\sum_{m=l+1}^n k_m^2\right) \right),
\end{align*}
which is negative since $k_1(p)=\cdots=k_l(p)<0.$ Thus we derive a contradiction with \eqref{41}, and the conclusion of the theorem follows.\qed

\section{\bf Gap theorem for Mean Convex $\lambda$-hypersurfaces}\label{5}
In \cite{CW18}, Cheng and Wei proved a rigidity theorem for $\lambda$-hypersurfaces under some curvature assumptions. Their result is an application of the arguments that Huisken applied in \cite{H90} and \cite{H93}. (Note that the definition of $\lambda$-hypersurfaces in \cite{CW18} is different from that in this article by a constant. The sign convention of the second fundamental form in \cite{CW18} is also different from ours.) We use the maximum principle to give a proof of the theorem here following the ideas in \cite{H90}. In the mean convex case, we can use theorem \ref{main} to derive a gap theorem when $\lambda\le 0.$

\begin{thm}[\cite{CW18}]\label{CW}
	Let $M^n$ be a smooth, closed, and embedded $\lambda$-hypersurface in $\bb R^{n+1}.$ If $H-\lambda\ge 0$ and $\lambda(2(H-\lambda)\tr A^3 + |A|^2) \le 0,$ then $M$ is a round sphere.
\end{thm} 
\noindent{\bf(Proof.)} By the maximum principle, we have $H-\lambda>0.$ Using \eqref{Hij}, \eqref{Simons}, and the $\lambda$-hypersurface equation \eqref{lambda}, we can derive
\begin{equation*}\label{DH}
\Delta H 
= \frac 12 H + \frac 12 \nabla_{x^T} H - (H-\lambda) |A|^2
\end{equation*}
and
\begin{equation*}
\Delta |A|^2
= 2|\nabla A|^2 + |A|^2 - 2|A|^4 + \frac 12\nabla_{x^T} |A|^2 - 2\lambda \tr A^3.
\end{equation*}
As a result,
\begin{align*}
\Delta\left(\frac{|A|^2}{(H-\lambda)^2} \right)
& = \frac{\Delta|A|^2}{(H-\lambda)^2}
- \frac{2|A|^2}{(H-\lambda)^3}\Delta H
- \frac{4}{(H-\lambda)^3}\pair{\nabla |A|^2,\nabla H}
+ \frac{6|A|^2}{(H-\lambda)^4}|\nabla H|^2\\
& = \frac 1{(H-\lambda)^4}\left( 
2(H-\lambda)^2|\nabla A|^2 
+ \frac 12 H^2\nabla_{x^T}|A|^2
- H|A|^2\nabla_{x^T}H
\right)\\
& +\frac 1{(H-\lambda)^4}\left(
-\lambda(H-\lambda)\left( 2(H-\lambda)\tr A^3 + |A|^2  \right)
- 4(H-\lambda)\pair{\nabla |A|^2,\nabla H}
+ 6|A|^2 |\nabla H|^2
\right).
\end{align*}
Plugging in 
$$
|a_{ij}\nabla_l H - (H-\lambda)\nabla_l a_{ij}|^2
= |A|^2 |\nabla H|^2
+ |\nabla A|^2 (H-\lambda)^2
- (H-\lambda)\pair{\nabla H,\nabla |A|^2}
$$
and
$$
\nabla\left(\frac{|A|^2}{(H-\lambda)^2} \right)
= \frac{\nabla |A|^2}{(H-\lambda)^2} - \frac{2|A|^2}{(H-\lambda)^3}\nabla H,
$$
we finally obtain
\begin{align*}
\Delta\left(\frac{|A|^2}{(H-\lambda)^2} \right)
& = \frac{2}{(H-\lambda)^2}\left(
|a_{ij}\nabla_l H - (H-\lambda)\nabla_l a_{ij}|^2
-\frac 12\lambda(H-\lambda)\left( 2(H-\lambda)\tr A^3 + |A|^2  \right)
\right)\\
& + \pair{- \frac{2}{H-\lambda}\nabla H
+ \frac {x^T}2,
\nabla\left(\frac{|A|^2}{(H-\lambda)^2} \right)}.
\end{align*}
By our assumptions, we have
$$|a_{ij}\nabla_l H - (H-\lambda)\nabla_l a_{ij}|^2
-\frac 12\lambda(H-\lambda)\left( 2(H-\lambda)\tr A^3 + |A|^2  \right)\ge 0,$$
so the maximum principle implies $|A|^2 = C(H-\lambda)^2$ for some constant $C$ and that
$$|a_{ij}\nabla_l H - (H-\lambda)\nabla_l a_{ij}|^2
-\frac 12\lambda(H-\lambda)\left( 2(H-\lambda)\tr A^3 + |A|^2  \right)= 0.$$
In particular, we have
$$|a_{ij}\nabla_l H - (H-\lambda)\nabla_l a_{ij}|^2=0.$$
This tells us that the anti-symmetric part of this tensor also vanishes, which implies
\begin{equation}\label{aijajl}
|a_{ij}\nabla_l H -a_{il}\nabla _j H|^2=0
\end{equation}
by the Codazzi equation. 

Now we assume $M$ is not a round sphere. Then we can find a point $p\in M$ at which $\nabla H\neq 0.$ If we take a local frame $e_1,\cdots,e_n$ such that $e_1=\frac{\nabla H}{|\nabla H|}$ at $p,$ then \eqref{aijajl} implies
$$|\nabla H|^2\left(|A|^2-\sum_{i=1}^n a_{1i}^2 \right)=0$$
at $p.$ Since $\nabla H(p)\neq 0,$ we get $|A|^2-\sum_{i=1}^n a_{1i}^2=0$ at $p.$ As a result,
$$
\sum_{i=1}^n a_{1i}^2
= |A|^2
= \sum_{j=1}^n\sum_{k=1}^n a_{jk}^2,
$$
which implies $a_{jk}=0$ if $(j,k)\neq (1,1).$ Thus we get $|A|^2=a_{11}^2=H^2.$ This along with the fact that $|A|/(H-\lambda)$ is constant implies that $H$ is constant, which leads to a contradiction since we assume $M$ is not a round sphere.\qed

We remark that when $\lambda=0,$ the calculations above reduce to those in \cite{H90}. Now we can use theorem \ref{CW} to prove our gap theorem \ref{gap} for mean convex $\lambda$-hypersurfaces. We state theorem \ref{gap} here again.
\begin{thm}
	Let $M^n$ be a smooth, closed, and embedded $\lambda$-hypersurface in $\bb R^{n+1}.$ If $\lambda\le 0$ and the mean curvature of $M$ satisfies
	$$0\le H\le \frac{\sqrt{\lambda^2+2}+\lambda}{2},$$
	then $M$ is a round sphere.
\end{thm}
\noindent{\bf (Proof.)} By theorem \ref{main}, we know that $M$ is convex. That is, $k_i\ge 0$ for all $i=1,\cdots,n.$ In particular, this implies
$$
H|A|^2 = \left(\sum_{i=1}^n k_i \right)\cdot \left(\sum_{j=1}^n k_j^2 \right)
\ge \sum_{m=1}^n k_m^3 
= -\tr A^3.
$$
On the other hand, by the upper bound of $H,$ we can conclude that $1\ge 2(H-\lambda)H.$ Combining these gives 
\begin{align*}
|A|^2 \ge 2(H-\lambda)H|A|^2
\ge -2(H-\lambda)\tr A^3,
\end{align*}
which implies $\lambda(2(H-\lambda)\tr A^3 + |A|^2) \le 0.$ Applying theorem \ref{CW}, the conclusion follows.\qed

\end{document}